%% file: root.tex
\documentclass[letterpaper, 10 pt, conference]{ieeeconf}
\usepackage{comment}
\IEEEoverridecommandlockouts                              % This command is only needed if 
\usepackage{graphics} % for pdf, bitmapped graphics files
\usepackage{epsfig} % for postscript graphics files
\usepackage{times} % assumes new font selection design installed
\usepackage{amsmath} % assumes amsmath package installed
\usepackage{amssymb}  % assumes amsmath package installed
\usepackage{xcolor} % add text colors
\usepackage{xfrac} % diagonal fraction
\usepackage{pgffor}
\usepackage{subcaption}
\usepackage{booktabs}
\usepackage{cite}
\usepackage{breakcites}

\usepackage{pgfplots}
  \pgfplotsset{compat=newest}
  \usetikzlibrary{plotmarks}
  \usetikzlibrary{arrows.meta}
  \usepgfplotslibrary{patchplots}
  \usepackage{grffile}

\newtheorem{proposition}{Proposition}

\title{\LARGE \bf
Balancing Safety and Traffic Throughput \\ in Cooperative Vehicle Platooning
}

\author{Stanley W. Smith$^1$, Yeojun Kim$^2$, Jacopo Guanetti$^2$, \\ Alexander A. Kurzhanskiy$^3$, Murat Arcak$^1$, and Francesco Borrelli$^2$ % <-this % stops a space
\thanks{This work was supported in part by an NDSEG Graduate Fellowship and NSF grant CNS-1545116.} % <-this % stops a space
\thanks{$^{1}$Dept. of Electrical Engineering and Computer Sciences, University of California, Berkeley
        {\tt\footnotesize $\{$swsmth,arcak$\}$@eecs.berkeley.edu.}}%
\thanks{$^{2}$Dept. of Mechanical Engineering, University of California, Berkeley
        {\tt\footnotesize $\{$yk4938, jacopoguanetti, fborrelli$\}$@berkeley.edu.}}%
\thanks{$^{3}$California Partners for Advanced Transportation Technology, University of California, Berkeley
        {\tt\footnotesize akurzhan@berkeley.edu.}}%
}

\begin{document}

\maketitle
\thispagestyle{empty}
\pagestyle{empty}

%%%%%%%%%%%%%%%%%%%%%%%%%%%%%%%%%%%%%%%%%%%%%%%%%%%%%%%%%%%%%%%%%%%%%%%%%%%%%%%%
\begin{abstract}
In this paper we propose  a distributed model predictive control architecture to coordinate the longitudinal motion of a vehicle platoon at a signalized intersection. Our control approach is cooperative; we use vehicle-to-vehicle (V2V) communication in order to maintain small inter-vehicle distances and correspondingly achieve large vehicle throughput at an intersection. We study the trade-off between safety and road throughput for this problem. In particular, we present the link between traffic efficiency gains in terms of throughput, and  safety of the connected platoon measured as trust on the predicted speed profile of other vehicles in the network.
\end{abstract}

%%%%%%%%%%%%%%%%%%%%%%%%%%%%%%%%%%%%%%%%%%%%%%%%%%%%%%%%%%%%%%%%%%%%%%%%%%%%%%%%
\section{Introduction}

Research on connected and automated vehicles involves perception, communication, control and decision making in a variety of driving scenarios.
The most widely adopted technology is Adaptive Cruise Control (ACC), a driving assistance system that maintains a cruising speed set by the driver when the lane is free, and reduces the speed to maintain a safe distance when a vehicle in front is detected.
ACC systems typically use radar and camera for perception of the front vehicle.
Cooperative ACC (CACC) enhances ACC with V2V communication, which has the potential to greatly mitigate traffic congestion, prevent accidents, and reduce reliance on traffic lights \cite{V2Varticle}. 

A prominent technology for V2V communication is Dedicated Short Range Communication (DSRC) \cite{kenney2011dedicated}, which enables very low latency communication between vehicles.
In the simplest case, V2V communication can be used to share the current state of each vehicle, including speed, acceleration, heading; this is beneficial to overcome some limitations of perception systems (noise, delays, and obstacles out of sight) and to react promptly to changes in the front vehicle(s) behavior (for instance, to changes in acceleration). V2V can also be used for more advanced forms of coordination, by sharing forecasts of each vehicle's behavior and by establishing vehicle formations in real-time.

Early research on connected vehicle technology studied \emph{platoons}, strings of vehicles that travel in a coordinated manner, generally at the same speed and at short distance.
Coordination on busy or congested highway segments has a large potential for enhancing the vehicle throughput, reducing congestion, and preventing accidents.
In 1994 and 1997, the California PATH team demonstrated a platoon of vehicles driving on the I-15 in San Diego, CA \cite{shladover2007path}. The vehicles were equipped with radars and also used V2V communication to enable cooperative driving \cite{chang1991experimentation}.  More recent demonstrations, such as \cite{milanes2014cooperative, ploeg2012introduction, ploeg2016introduction},
%we should add more refs here that Milanes et al. cite
aimed to advance platooning to urban and realistic traffic situations.

A less explored benefit of vehicle platooning is its ability to dramatically increase vehicle throughput at intersections.
Indeed, the throughput, measured as the number of vehicles that pass through a given road segment (e.g, an intersection) per hour (vph), can potentially be \emph{doubled} by forming platoons \cite{lioris2017platoons}.
To fully realize this potential, small inter-vehicle distances must be maintained within the platoons passing through an intersection.
However, doing so while accelerating or decelerating constitutes a challenging motion coordination problem with trade-offs between safety and performance.

In this paper, we propose a distributed Model Predictive Control (MPC) approach for the coordination of a vehicle platoon.
We investigate the performance of our approach via vehicle throughput, and then demonstrate the trade-off between throughput and safety (conservativeness of the controller), where the latter is quantified by a parameter in our MPC formulation.

% Literature on autonomous intersections ...
% Open problem: trade off between safety and throughput, when crossing a signalized intersection

% - Our contributions: \\
% - obtain baseline for vehicle throughput with a conservative control design guaranteeing safety (only trust radar) \\
% - improve throughput from baseline by (some possibilities): using accel. info from preceding vehicle, info from leader (PFL topology), allowing small safety violations, making some assumptions, etc

\section{Platoon Model}

\subsection{Vehicle dynamics}
% - 1 leader vehicle, 3 follower vehicles \\
% - vehicle model (homogeneous across platoon) \\
% - discretized dynamics for MPC \\
In this paper we consider a platoon comprised of $N$ vehicles. The vehicle positioned at the front of the platoon is referred to as the leader vehicle, the subsequent vehicles behind it as follower vehicles $1, \dots, N-1$.

We model the longitudinal dynamics \cite{guzzella2007vehicle} of the leader vehicle as
\begin{subequations} \label{leaderDyn}
\begin{align}
\dot{p}^L &= v^L, \\
\dot{v}^L &= \frac{1}{M} \left(\frac{T_w^L}{R_w} - F_f^L \right), \label{accel1}
\end{align}
\end{subequations}
where the states $p^L$ and $v^L$ are the absolute position and velocity of the leader vehicle, respectively, and the input $T_w^L$ is the wheel torque. $M$ is the mass of the vehicle, $R_w$ is the wheel radius, and $F_f^L$ is the frictional force, given by
\begin{equation} \label{frictionalForce}
F_f^L = M g (sin(\theta) + c_r cos(\theta)) + \frac{1}{2} \rho A c_x (v^L)^2 
\end{equation}
where $g$ is the gravitational constant, $\theta$ is the road grade, $c_r$ is a rolling coefficient, $\rho$ is air density, $A$ is the area of the vehicle, and $c_x$ is an air drag coefficient. In this paper we do not consider road grade, and thus $\theta = 0$ for all $t \geq 0$.

While it is reasonable to assume the velocity of the leader $v^L$ is accessible, to determine the position of the leader $p^L$ one must perform localization using sensor measurements (e.g. differential GPS), which is difficult in general. Once $p^L$ has been estimated, the states in \eqref{leaderDyn} can be used to compute a control input for $T_w^L$, or to communicate state information to nearby vehicles. We write the dynamics \eqref{leaderDyn} for the leader vehicle concisely as
\begin{equation} \label{leaderConcise}
\dot{x}^L = f^L(x^L, u^L)
\end{equation}
where $x^L := [p^L;\ v^L]$ and $u^L := T^L_w$.

\begin{table}[t]
    \caption{Model Parameters}
    \label{modelParameters}
    \centering
    \begin{tabular}{c l l l}
    \toprule
         $M$ & vehicle mass & kg & 1722 \\
         $A$ & vehicle reference area & $\text{m}^2$ & 2.6292 \\
         $\rho$ & air density & $\text{kg}/\text{m}^3$ & 1.206 \\
         $c_x$ & vehicle drag coefficient & - & 0.2047 \\
         $c_r$ & vehicle roll coefficient & - & 0.0106 \\
         $\Delta t$ & sampling time & s & 0.1 \\
    \bottomrule
    \end{tabular}
\end{table}

Each follower vehicle is assumed to have a radar which accurately measures the distance $h^i$ (for follower vehicle $i = 1, \dots, N-1$) between itself and its preceding vehicle, as well as their relative speed. The distance $s^i$ between each follower vehicle and the leader vehicle is also modelled (it is estimated using V2V messages, discussed further in Section \ref{communication}). The longitudinal dynamics of the follower vehicles are modelled as
\begin{subequations} \label{followerDyn}
\begin{align}
\dot{p}^i &= v^i, \\
\dot{s}^i &= v^L - v^i, \\
\dot{h}^i &= v^{i-1} - v^i, \\
\dot{v}^i &= \frac{1}{M} \left( \frac{T_w^i}{R_w} - F_f^i \right), \quad i = 1, \dots, N-1, \label{accel2}
\end{align}
\end{subequations}
where we let $v^0 = v^L$ so the relative position of follower vehicle $1$ is given with respect to the leader vehicle (note also that $s^1$ = $h^1$). The frictional force $F_f^i$ for follower vehicle $i$ is also modelled as in \eqref{frictionalForce}, with $v^i$ replacing $v^L$. Vehicle dynamics are assumed to be homogeneous within the platoon so that model parameters do not vary between vehicles (see Table \ref{modelParameters}). We similarly write \eqref{followerDyn} as
\begin{align} \label{followerConcise}
\dot{x}^i = f(x^i, u^i, w^i)
\end{align}
where $x^i := [p^i;\ s^i;\ h^i;\ v^i]$, $u^i := T^i_{w}$, and $w^i := [v^L;\ v^{i-1}]$. Note that $w^i$, containing the velocities of the leader and preceding vehicles, appears as a disturbance here. In our MPC formulation, planned velocity trajectories of the leader/preceding vehicles (received via V2V communication) are used as disturbance previews.

Next, we wish to obtain linear discrete-time models of \eqref{leaderConcise} and \eqref{followerConcise}. For the follower vehicles we first linearize \eqref{followerConcise} about the nominal velocity $v^i = v^i_0$, resulting in
\begin{equation} \label{linearized}
\dot{x}^i = \overline{A} x^i + \overline{B} u^i + \overline{E} w^i
\end{equation}
where $\overline{A} = \left. \frac{\partial f}{\partial x^i} \right|_{v^i = v^i_0}$, $\overline{B} = \frac{\partial f}{\partial u^i} $, and $\overline{E} = \frac{\partial f}{\partial w^i}$, and the matrix $\overline{A}$ is a function of the velocity $v^i_0$ due to the velocity squared term in \eqref{frictionalForce}. We then discretize \eqref{linearized} with time step $\Delta t = 0.1s$, and obtain
\begin{equation} \label{linearDisDyn}
x^i (k+1) = A x^i(k) + B u^i(k) + E w^i(k)
\end{equation}
where $A = e^{\overline{A} \Delta t}$, $B = \int_0^{\Delta t} e^{\overline{A} \tau} \overline{B} d\tau$, and $E = \int_0^{\Delta t} e^{\overline{A} \tau} \overline{E} d\tau$. Since $\overline{A}$ is a function of $v^i_0$, computing $A$, $B$, and $E$ requires integrating symbolic expressions - the resulting matrices are stored as symbolic functions of $v^i_0$. At each time step during simulation, we substitute the current velocity $v^i_0$ into these functions to obtain the appropriate model to be used for MPC. The same procedure is used to obtain the model of the leader vehicle dynamics
\begin{equation}
x^L (k+1) = A^L x^L(k) + B^L u^L(k)
\end{equation}
which is also parameterized by velocity.

\begin{figure}
\centering
\includegraphics[width = \columnwidth]{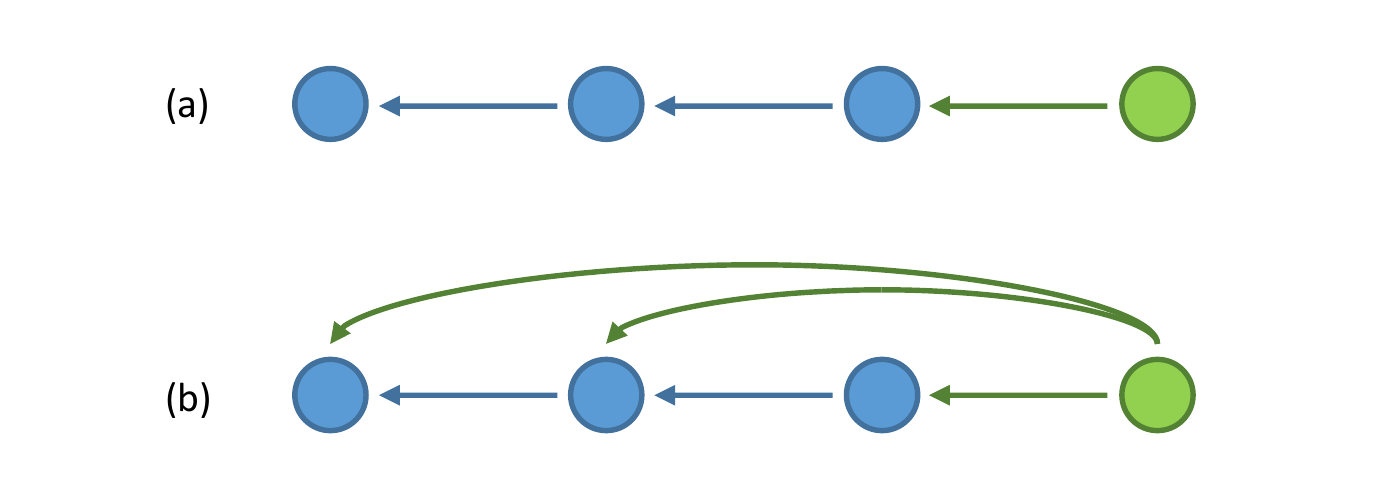}
\caption{Two information flow topologies for $N = 4$: the predecessor following (top) and predecessor following leader (bottom). Messages travel in the direction of the arrows.}
\end{figure}

\subsection{Vehicle-to-vehicle communication} \label{communication}
% - information flow topologies (predecessor following $\&$ predecessor following leader) \\
% - message content \\
% - modeling/accounting for delays (briefly): delay effects on safety?
To permit cooperation, we assume each vehicle is equipped with vehicle-to-vehicle (V2V) communication capabilities. Many information flow topologies have been proposed for the coordination of vehicle platoons (see, e.g., \cite{guanetti2018control}). In the predecessor following topology (Figure 1a), each follower vehicle receives a message only from the preceding vehicle. The predecessor following leader topology (Figure 1b) also includes communication arcs between the leader vehicle and followers 2 and 3, which can enable improvements in performance. Note that here the arcs are unidirectional, so that information only flows backwards away from the leader of the platoon. Messages are sent at a rate of ten messages per second - the same time step used for MPC. In this paper, we use the predecessor following leader topology.

Messages sent at time $t_{sent}$ contain the following information
\begin{align}
m^L &= [t_{sent};\ p^L(t|t);\ v^L(t|t);\ \dots \ v^L(t+N_p|t) ] \nonumber \\
m^i &= [t_{sent};\ v^i(t|t);\ \dots \;\ v^i(t+N_p|t) ] \label{V2Vmessages}
\end{align}
where $m^L$ and $m^i$ are the messages sent by the leader and follower vehicle $i$, respectively, and $t_{sent}$ is the time stamp contained in the message. Furthermore, $t$ is the current time step, and $v^L(k|t)$ is the planned velocity of the leader at time step $k$ (obtained by solving an MPC problem at time step $t$). For example, $p^L(t|t)$ refers to the current position of the leader. The planning horizon of the MPC is $N_p$ time steps, so that $v^L(t+N_p|t)$ is the terminal velocity of the leader. The notation is the same for the follower vehicles. Delays in V2V communication are also modelled: a vehicle receiving a delayed message can over-estimate the number of time steps of the delay via
\begin{equation}
d = \text{ceil}\left( \frac{t_{received} - t_{sent}}{\Delta t} \right) \label{delayComputation}
\end{equation}

\section{Traffic Throughput}
\subsection{Problem Statement} \label{throughputProblem}
In this paper, we consider a scenario where the platoon accelerates from a stop in response to a traffic light cycling from red to green. At the scenario start, the vehicle states are
\begin{align} \label{initialState}
x^L(0) &= [0; \ 0], \nonumber \\
x^i(0) &= [-s \cdot i; \ s \cdot i; \ s; \ 0], \quad i = 1, \dots, N-1,
\end{align}
where $p^L(0) = 0$ is assumed to be the position of the intersection stop bar and $s \in \mathbb{R}_{>0}$ is the initial distance between vehicles in the platoon.

To measure the performance of our control design in terms of vehicle throughput, we choose a point $\ell$ in the intersection, and define $t_L$ and $t_{N-1}$ to be the smallest time instants in seconds such that $p^L(t_L) \geq \ell$ and $p^{N-1}(t_{N-1}) \geq \ell$, respectively. We note that $\ell$ should be chosen sufficiently large to avoid taking measurements of $t_L$ and $t_{N-1}$ while the platoon is still in transient behavior. Then, the throughput in vehicles per hour can be estimated as
\begin{equation} \label{throughput}
\text{throughput} \approx 3600 \cdot \left( \frac{N-1}{t_{N-1} - t_L} \right).
\end{equation}
The goal is to achieve a high level of throughput while keeping the platoon safe. In the MPC formulation for the follower vehicles (Section \ref{followerMPCsection}), a large throughput is assured by penalizing deviations from a small desired distance between the vehicles in the platoon. % Of course, it is also possible to clear an intersection quickly if all vehicles apply a large acceleration immediately; doing so, however, creates an uncomfortable and potentially unsafe experience for the passengers. We avoid this behavior by including penalties on vehicle jerk and also terminal set constraints for safety in our MPC formulation.

\section{MPC Formulation}

\subsection{Leader MPC}
% - MPC only trusting radar gives performance baseline \\
% - throughput also limited by acceleration capabilities \\
% - for clarity, we add hats to variables that are estimated \\
% - theorem/proposition at the end that states safety is guaranteed

We now present our control design for the platoon. We begin with a control design which guarantees safety of all vehicles in the platoon, but which will not necessarily result in a high level of throughput at an intersection - providing a baseline for comparison. We then show how to improve upon the baseline by slightly relaxing the safety constraints, and then further explore the trade-off between traffic throughput and safety.

Throughout the paper we use a distributed MPC approach to coordinate the motion of the platoon. In particular, each vehicle obtains its control input at each time step by solving its own finite horizon MPC problem. The control objective of the leader vehicle is to attain and preserve a desired velocity $v_{des}$ - its MPC problem is given by
\begin{subequations}
\begin{align}
& \underset{u(\cdot|t)}{\text{min}}
& J^L & = (v^L(t+N_p|t) - v^L_{des})^2 \label{leaderObj1} \\
& & & + \alpha \sum_{k = t}^{t + N_p - 1} ( u^L(k+1|t) - u^L(k|t) )^2 \label{leaderObj2} \\
& \text{s.t.} & & x^L(k+1|t) = A^L x^L(k|t) + B^L u^L(k|t), \label{leaderConstr1} \\
& & & v_{min} \leq v^L(k|t) \leq v_{max}, \label{leaderConstr2} \\
& & & u_{min} \leq u^L(k|t) \leq u_{max}, \label{leaderConstr3} \\
& & & x^L(t|t) = x^L(t), \label{leaderConstr5} \\
& & & \forall k = t, \dots, t+N_p-1, \nonumber
\end{align}
\end{subequations}
where $N_p$ is the horizon of the controller in time steps, and $x^L(k|t)$ and $u^L(k|t)$ are the planned state and input of the leader vehicle at time step $k$, computed at time step $t$, respectively. Note that, for simplicity, this MPC problem assumes the leader vehicle has no obstacles. In practice, this formulation can be easily augmented to incorporate obstacles or a lead vehicle (e.g., as ACC does).

The leader objective function $J^L$ includes two major components: \eqref{leaderObj1} penalizes deviations from the desired velocity $v_{des}$ at the end of the planning horizon and \eqref{leaderObj2} penalizes vehicle jerk. The parameter $\alpha > 0$ appearing in $J^L$ can be increased, placing a higher cost on vehicle jerk, if smoother acceleration profiles are desired. Furthermore, the constraints for the leader vehicle MPC problem are the state dynamics \eqref{leaderConstr1}, upper and lower bounds on the velocity and wheel torque \eqref{leaderConstr2} $\&$ \eqref{leaderConstr3}, and the initial condition \eqref{leaderConstr5}. 
%We note that the slack variable $\delta \geq 0$ allows violation of the slew rate constraint \eqref{leaderConstr4} if needed.

\subsection{Follower MPC} \label{followerMPCsection}
For the follower vehicles, the primary objective is to maintain a desired distance $s^i_{des}$ to the leader vehicle at all times, while also ensuring safety if the preceding vehicle decelerates at the maximum rate. For each follower vehicle $i$, the corresponding MPC problem is given by 
\begin{subequations} \label{followerMPC}
\begin{align}
& \underset{u(\cdot|t)}{\text{min}}
& J^i & = \sum_{i=t}^{t+N_p} (s^i(k|t) - s^i_{des})^2 \label{followerObj1} \\
& & & + \alpha \sum_{i=t}^{t+N_p-1} (u^i(k+1|t) - u^i(k|t))^2 \label{followerObj2} \\
& \text{s.t.} & & x^i(k+1|t) = \label{followerConstr1} \\
& & & \quad A^i x^i(k|t) + B^i u^i(k|t) + E^i \hat{w}^i(k), \nonumber \\
& & & v_{min} \leq v^i(k|t) \leq v_{max}, \label{followerConstr2} \\
& & & h_{min} \leq h^i(k|t), \label{minDistanceMPC} \\
& & & u_{min} \leq u^i(k|t) \leq u_{max}, \label{followerConstr3} \\
& & & x^i(t|t) = \hat{x}^i(t), \label{followerConstr5} \\
& & & \forall k = t, \dots, t+N_p-1, \nonumber \\
& & & \begin{bmatrix}
h^i(t+F|t) \\ v^i(t+F|t)
\end{bmatrix} \in C(\hat{v}^{i-1}(t+F)). \label{followerConstr6}
\end{align}
\end{subequations}

\begin{table}
    \caption{MPC Parameters}
    \label{modelParameters}
    \centering
    \begin{tabular}{c l l l}
    \toprule
         $h_{des}$ & desired distance & m & 9 \\
         $h_{min}$ & minimum distance & m & 6.5 \\
         $v_{min}$ & minimum velocity & m$/$s & 0 \\
         $v_{max}$ & maximum velocity & m$/$s & 30 \\
         $v_{des}$ & desired velocity & m$/$s & 15.64 \\
         $u_{max}$ & maximum wheel torque & Nm & 1500 \\
         $u_{min}$ & minimum wheel torque & Nm & -2000 \\
         $\Delta u_{max}$ & maximum slew rate & Nm$/$s & 250 \\
         $N_p$ & MPC horizon & - & 20 \\
    \bottomrule
    \end{tabular}
\end{table}

The follower vehicle objective function $J^i$ includes a term \eqref{followerObj1} penalizing deviations from the desired distance to the leader, as well as penalties on jerk and the slack variable $\delta$, similar to $J^L$. Here, the desired distance to the leader vehicle $s^i_{des}$ is defined as
\begin{equation}
s^i_{des} := h_{des} \cdot i
\end{equation}
where $h_{des}$, the desired distance to the preceding vehicle, is a design parameter.

Constraints \eqref{followerConstr1} - \eqref{followerConstr5} are similar to those in the leader MPC problem - except here we also include a minimum distance constraint \eqref{minDistanceMPC} for safety. The set constraint \eqref{followerConstr6} is discussed in Section \ref{safetyVsThroughput}. The disturbance preview $\hat{w}^i(k)$ and the initial state $\hat{x}^i(t)$ are estimated using received V2V messages. Here, $\hat{w}^i(k) = [\hat{v}^L(k);\ \hat{v}^{i-1}(k)]$ contains estimates for the velocity of the leader/preceding vehicle at time step $k$. Since the planned trajectory of each vehicle changes at every time step, we estimate future velocities using the most recent V2V message. For example, for the leader vehicle we set
\begin{align}
&\hat{v}^L(k) = \label{velocityEstimate} \\
& \ \begin{cases}
v(k|t-d), & k = t, \dots, t-d+N_p, \\
v(t-d+N_p|t-d), & k = t-d+N_p+1, \dots, t+N_p,
\end{cases} \nonumber
\end{align}
where the delay $d$ is computed as in \eqref{delayComputation}. In other words, the MPC problem \eqref{followerMPC} assumes the leader maintains a constant velocity beyond its planned trajectory. The estimate for the velocity of the preceding vehicle is set differently, and is discussed in Section \ref{safetyVsThroughput}. The initial state $\hat{x}^i(t)$ contains the current distance to the leader vehicle, estimated using the most recent message from the leader as
\begin{align}
& \hat{s}^i(t) = \label{positionEstimate} \\
& \quad \left( p^L(t-d|t-d) + \Delta t \sum_{k=0}^{d-1} v(t-d+k|t-d) \right) - p^i(t) \nonumber
\end{align}
where the term in parentheses approximates $p^L(t)$. We include delay in \eqref{velocityEstimate} and \eqref{positionEstimate} to make these estimates robust to communication latencies. The leader could also communicate its entire planned position trajectory to eliminate the need for approximation - however, we do not do this so the leader does not need to send large amounts of data. The remaining states in $\hat{x}^i(t)$ are assumed to be available from sensor data.

\textbf{\begin{figure}[t]
\centering
\input{tikz/terminal_set.tex}
\caption{Set $\mathcal{C}(v_0)$ with $v_0 = 7.5 \text{m}/\text{s}$, $a_{min} = -3.218 \text{ m}/\text{s}^2$, and $h_{min} = 6.5$ m.}
\label{terminalSet}
\end{figure}
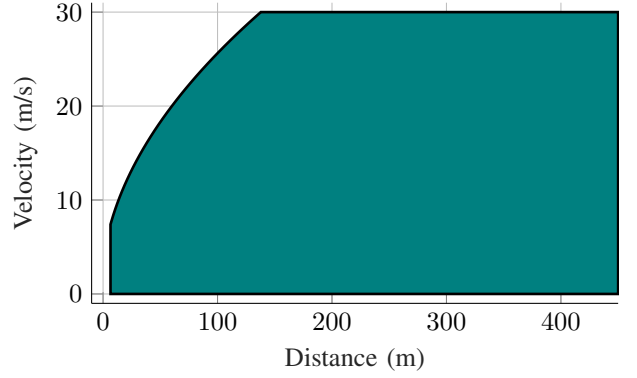}

\begin{figure*}
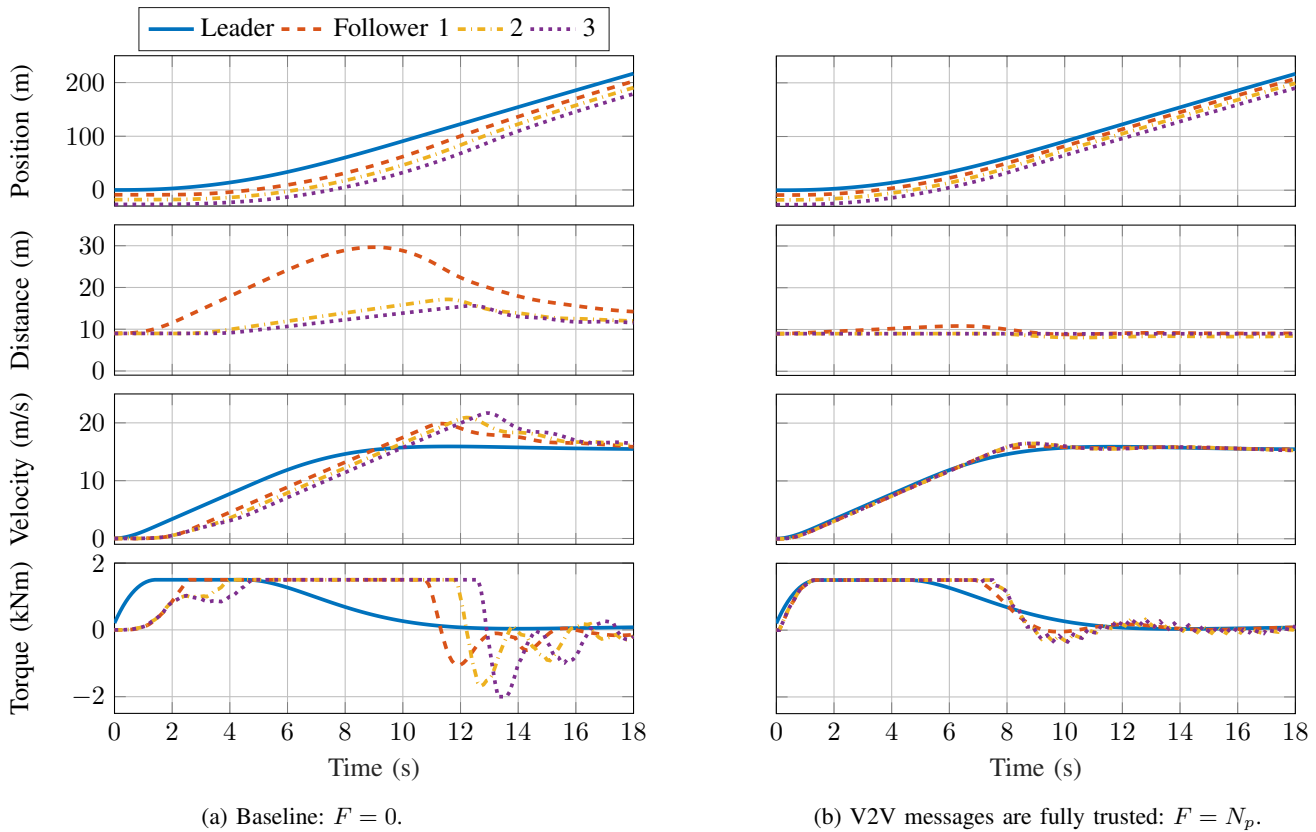

    \centering
    \begin{subfigure}{0.45\textwidth}
        \centering
        \input{tikz/CACC_no_trust_new2.tex}
        \caption{Baseline: $F = 0$.}
        \label{noTrust}
    \end{subfigure}
    \hfill
    \begin{subfigure}{0.45\textwidth}
        \centering
        \vspace{0.6475cm}
        \input{tikz/CACC_full_trust_new2.tex}
        \caption{V2V messages are fully trusted: $F = N_p$.}
        \label{fullTrust}
    \end{subfigure}
    \caption{Full simulation results (with a platoon of size $N = 4$) for the baseline (left), where only radar data is trusted, and also for where the V2V messages are fully trusted (right).}
    \label{fullSimulation}
\end{figure*}

\subsection{Safety vs. Throughput} \label{safetyVsThroughput}
The set constraint \eqref{followerConstr6} ensures each follower can persistently avoid violation of a minimum distance constraint
\begin{equation} \label{minDistance}
h^i(k) \geq h_{min}, \quad i = 1, \dots, N-1
\end{equation}
and velocity constraint \eqref{followerConstr2} if the preceding vehicle decelerates. To compute the set numerically, we approximate the dynamics of each vehicle with the following kinematic model
\begin{align}
p(k+1) &\approx p(k) + v(k) \Delta t + \frac{1}{2} a(k) \Delta t^2, \label{kinematicDyn} \\
v(k+1) &\approx v(k) + a(k) \Delta t, \nonumber
\end{align}
where $a(k)$ is acceleration. Suppose, starting at the current time step, the preceding vehicle decelerates from the velocity $v_0$ with $a(k) = a_{min} < 0$, until coming to a stop in $k_s$ time steps. Then, linearity of \eqref{kinematicDyn} enables us to use the Multi-Parametric Toolbox \cite{herceg2013multi} to compute the set of safe states, defined as
\begin{align}
\mathcal{C}(v_0) := \nonumber
\{\ &[h^i(0);\ v^i(0)] \mid \exists\ u^i(k) \text{ for } k \geq 0 \quad \text{s.t} \nonumber \\
& \text{vehicle dynamics as in \eqref{kinematicDyn} }, \nonumber \\
& v^{i-1}(0) = \tilde{v}_0, \nonumber \\
& a^{i-1}(k) = a_{min}, \quad \forall k = 0, \dots, k_s - 1, \nonumber \\
& v^{i-1}(k_s) = 0, \nonumber \\
& h^i(k) \geq h_{min}, \quad \forall k \geq 0 \}. \label{safeSet}
\end{align}
We note that in addition to $v_0$, the set $\mathcal{C}$ also depends on $a_{min}$ and $h_{min}$. To under-approximate the number of time steps it will take the preceding vehicle to come to a stop, we take
\begin{equation}
k_s = \text{floor} \left( \frac{v_0}{|a_{min}| \cdot \Delta t} \right)
\end{equation}
where $\text{floor}(\cdot) : \mathbb{R} \to \mathbb{Z}$ is the floor function. The velocity $\tilde{v}_0$ in \eqref{safeSet} is a slight under-approximation of $v_0$; starting at velocity $\tilde{v}_0$, the preceding vehicle can stop in \textit{exactly} $k_s$ time steps by applying the maximum deceleration $a_{min}$. Hence, we have $\tilde{v}_0 = |a_{min}| \cdot \Delta t \cdot k_s$. We compute a collection of sets offline corresponding to values of $v_0 \in [v_{min}, v_{max}]$, and then select the proper set to be used in our MPC problem at each time step during simulation. The full details on computing the sets are available in \cite{lefevre2016learning}. It is mainly important to note that the set depends on $v_0$, the velocity of the preceding vehicle before braking (Figure \ref{terminalSet} shows $\mathcal{C}(v_0)$ for $v_0 = 7.5$ m/s).

During simulation, the preceding vehicle is assumed to decelerate at rate $a_{min}$ starting at time $t+F$. We use $v_0 = \hat{v}^{i-1}(t+F)$ to select the set for \eqref{followerConstr6}, therefore the constraint \eqref{followerConstr6} is imposed. We refer to $F$ as our \textit{trust horizon}, which indicates the number of time steps of the preceding vehicle's planned velocity trajectory that are trusted. The velocity estimate at time $t+F$ is obtained using radar ($F = 0$), or the most recent V2V message from the preceding vehicle ($F > 0$). We note that if F = 0, we must place the constraint (14h) at time step t + 1, since the current state cannot be affected. We then set our estimated velocity trajectory for the preceding vehicle as
\begin{align} \label{brakingMessage}
&\hat{v}^{i-1}(k) = \\
& \resizebox{\linewidth}{!}{ \ $\begin{cases}
\text{as in } \eqref{velocityEstimate}, & k = t, \dots, t+F-1, \\
\tilde{v}_0, & k = t+F, \\
\max(0, \hat{v}^{i-1}(k-1) - a_{min} \cdot \Delta t), & k = t+F+1, \dots, t+N_p
\end{cases} $ } \nonumber
\end{align}
to incorporate the assumed braking behavior. We also use \eqref{brakingMessage} to set the leader vehicle's estimated velocity trajectory in the MPC problem for follower vehicle 1, since the leader vehicle of the platoon is its preceding vehicle.

Suppose we set $F = 0$ and $v_0 = \hat{v}^{i-1}(t)$ - i.e., we assume the preceding vehicle will begin braking immediately from its current velocity, measured via radar. Since this is the worst-case scenario, and we assume the current radar data is reliable, this guarantees safety of all vehicles in the platoon.

\begin{proposition}
Define the maximum deceleration as
\begin{align}
a_{min} := \frac{1}{M} \left( \frac{u_{min}}{R_w} - F_f^{max} \right) \label{maxDecel}
\end{align}
where
\begin{align}
& F_f^{max} := M g (sin(\theta) + c_r cos(\theta)) + \frac{1}{2} \rho A c_x (v_{max})^2. \label{maxFriction}
\end{align}
Since $a_{min}$ is a lower bound for the acceleration of the preceding vehicle (noting \eqref{accel1}, \eqref{accel2}), imposing the set constraint \eqref{followerConstr6} with $F = 0$, $v_0 = \hat{v}^{i-1}(t)$, and $a_{min}$ as in \eqref{maxDecel} ensures the MPC problem \eqref{followerMPC} is persistently feasible with respect to the minimum distance constraint \eqref{minDistance}. The proof is similar to that of Theorem 1 in \cite{lefevre2016learning} and is omitted.
\end{proposition}

Although setting $F = 0$ guarantees absolute safety of the platoon (i.e., \eqref{minDistance} holds), doing so results in poor traffic throughput. Thus, we use $F = 0$ as a baseline. To achieve higher throughput, we can extend the trust horizon by taking $F > 0$ and setting $v_0 = \hat{v}^{i-1}(t+F)$, which is estimated using the most recent V2V message. Doing so introduces some risk to the follower vehicles; however, we will see in Section \ref{simResults} that setting $F = 0$ is very restrictive and that much better throughput can be achieved by taking $F > 0$.
% \todo{(add theorem w/ upper bound on relative velocity if collision occurs..)}

% \subsection{Improving throughput}
% - make improvements to the baseline approach, possibly by slightly relaxing the conservative safety constraints, or using acceleration data received via V2V communication

\section{Simulation Results} \label{simResults}
We now present our simulation results for the signalized intersection scenario described in Section \ref{throughputProblem} for a platoon of size $N = 4$. For each simulation, the initial state of the vehicles is given by \eqref{initialState}, with the initial distance $s = 6.5$ m. The dynamics of each vehicle are modelled and simulated in MATLAB. A constant delay of $0.1$s, or exactly one time step, is also modelled for every communication arc in the platoon. The MPC problem for each vehicle is set up using Yalmip \cite{Lofberg2004} and solved using Gurobi \cite{gurobi}.

\begin{figure}
    \centering
    \input{tikz/tradeoff2.tex}
    \caption{Throughputs from \eqref{throughput} for various values of the trust horizon $F$.}
    \label{tradeoff}
\end{figure}
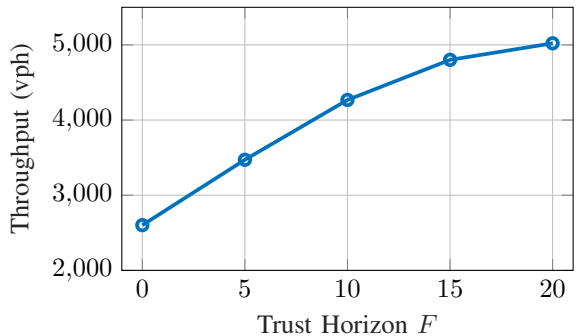

Full simulation results are shown in Figure \ref{fullSimulation} for the baseline, corresponding to $F = 0$, and also for where V2V messages are fully trusted, corresponding to $F = N_p$. In Figure \ref{noTrust}, we see that setting $F = 0$ results in large expansions in the inter-vehicle distances as the platoon accelerates. Furthermore, since only radar is utilized when $F = 0$, we also see some string instability. On the other hand, Figure \ref{fullTrust} shows that the inter-vehicle distances converge quickly to $h_{des} = 9$m when we set $F = N_p$, and that the instability is eliminated. Thus, we see that in order for platooning to enable a reasonable increase in traffic efficiency, each vehicle must trust a large portion of the planned trajectories received via V2V. However, doing so comes at the cost of exposing the vehicles in the platoon to some risk: potential slight violation of \eqref{minDistance} if a vehicle decelerates rapidly.

To further demonstrate the trade-off between safety and throughput, we vary the trust horizon parameter $F$ appearing in \eqref{safeSet}, and then estimate the corresponding throughput via \eqref{throughput}. Recall that $\ell$ is the position where we measure the crossing times of the leader and third follower vehicle to estimate throughput. Here, we choose $\ell = 30$m. Since the position of each vehicle is sampled every $0.1$s, we use linear interpolation to more accurately estimate both crossing times. As expected, Figure \ref{tradeoff} shows that the throughput increases as the trust horizon $F$ increases - the highest level of throughput is achieved when $F = 20$. % We hypothesize that this is because smoother vehicle trajectories are easier for the follower vehicles to react to, allowing the inter-vehicle distances to converge more quickly, and ultimately resulting in the highest level of throughput when $\alpha = 10^{-5}$ and $F = 20$.

% \subsection{Comparison}
% - comparison of the two approaches
% - copmarisons of conservativeness vs. throughput

\section{Conclusion}
In this paper we explored the trade-off between safety and throughput for vehicle platooning at an intersection. We formulated MPC problems to be solved by each vehicle in a distributed fashion in order for the platoon to: 1) achieve a desired velocity, 2) maintain small inter-vehicle distances, 3) ensure safety in the event that a vehicle decelerates rapidly. In particular, the MPC formulation included a parameter which quantifies the safety of the control design. In Section \ref{simResults} of the paper, we demonstrated the aforementioned trade-off by varying these parameters across multiple simulation runs, and then measuring the corresponding traffic throughput for each. Our results suggest that in order to achieve a reasonable increase in throughput, less restrictive safety constraints must be imposed on the vehicles in the platoon.

% In the future, it will be informative to quantify the risk associated with taking $F > 0$, measured via an upper bound on the relative velocity of the vehicles at the time of collision (if one occurs). Another research direction is to use an analytic form of the safe set $\mathcal{C}$, defined using the true dynamics of the vehicles. Since the constraints describing this set are nonlinear, a polyhedral inner-approximation could be used for the purpose of MPC.

\section*{Acknowledgements}
The authors would like to thank Pravin Varaiya and Roberto Horowitz for providing helpful suggestions.

\bibliographystyle{IEEEtran}

\end{document}

%% file: tikz/terminal_set.tex
\begin{tikzpicture}

\begin{axis}[%
width=7cm,
height=4cm,
at={(0cm,0cm)},
scale only axis,
xmin=-10,
xmax=450,
xlabel style={font=\color{white!15!black}},
xlabel={Distance (m)},
ymin=-1,
ymax=31,
ylabel style={font=\color{white!15!black}},
ylabel={Velocity (m/s)},
axis background/.style={fill=white},
axis x line*=bottom,
axis y line*=left,
xmajorgrids,
ymajorgrids,
legend style={legend cell align=left, align=left, draw=white!15!black}
]

\addplot[area legend, line width=1.0pt, draw=black, fill=teal]
table[row sep=crcr] {%
x	y\\
7.25623374835777	7.72323829109795\\
8.8652417265439	8.36684148239882\\
10.6029703431086	9.01044467370254\\
12.4694195952616	9.65404786409563\\
14.4645894879613	10.2976510553106\\
16.5884800182139	10.9412542462943\\
18.8410911863248	11.5848574371935\\
21.2224229931235	12.2284606282132\\
23.7324754384717	12.8720638193335\\
26.3712485192846	13.515667009785\\
29.1387422412918	14.1592702009501\\
32.0349565986398	14.8028733914987\\
35.0598915990849	15.4464765830961\\
38.2135472360496	16.0900797743072\\
41.4959235078422	16.7336829648511\\
44.9070204180571	17.3772861554658\\
48.4468379681966	18.0208893464181\\
52.1153761636851	18.6644925386051\\
55.9126349866019	19.308095728919\\
59.8386144501565	19.9516989196552\\
63.8933145501178	20.5953021101391\\
68.0767352874136	21.238905300521\\
72.3888766675227	21.8825084916134\\
76.8297386934937	22.5261116838057\\
81.3993213454924	23.1697148742218\\
86.0976246432801	23.8133180657018\\
90.9246485778967	24.4569212569785\\
95.8803931459806	25.1005244476246\\
100.964858360621	25.7441276393393\\
106.17804420058	26.3877308294533\\
111.519950693171	27.0313340213127\\
116.990577810164	27.6749372115394\\
122.589925575787	28.3185404029742\\
128.31799396564	28.9621435928461\\
134.174782998478	29.6057467832749\\
137.828870415269	30\\
% 19991.4883477995	30\\
% 19991.4883477995	0\\
450 30 \\
450 0 \\
6.5	-5.32907051820075e-15\\
6.5	7.40143669605203\\
}--cycle;
% \addlegendentry{data1}

\end{axis}
\end{tikzpicture}%

%% file: tikz/tradeoff2.tex
\definecolor{mycolor1}{rgb}{0.00000,0.44700,0.74100}%
\begin{tikzpicture}

\begin{axis}[%
width=6cm,
height=3.5cm,
at={(0cm,0cm)},
scale only axis,
xmin=-1,
xmax=21,
xlabel style={font=\color{white!15!black}},
xlabel={Trust Horizon $F$},
ymin=2000,
ymax=5500,
ylabel style={font=\color{white!15!black}},
ylabel={Throughput (vph)},
axis background/.style={fill=white},
xmajorgrids,
ymajorgrids,
legend style={legend cell align=left, align=left, draw=white!15!black}
]
\addplot [color=mycolor1, line width=1.4pt, mark=o, mark options={solid, mycolor1}]
  table[row sep=crcr]{%
0	2601.730059329\\
5	3471.63484898519\\
10	4267.71261714146\\
15	4801.1854804814\\
20	5021.17791622325\\
};
% \addlegendentry{data1}

\end{axis}
\end{tikzpicture}%